\documentclass[11pt,a4paper]{article}
\RequirePackage{amsmath}%
\RequirePackage{amsthm}%
\usepackage{amsfonts}%
\usepackage{amssymb}%
\usepackage{graphicx}%
\usepackage{cite}

\renewcommand{\Re}{\mathop{\mathrm{Re}}}
\renewcommand{\Im}{\mathop{\mathrm{Im}}}

\renewcommand{\i}{\mathrm{i}}

\renewcommand{\Re}{\mathop{\mathrm{Re}}}
\renewcommand{\Im}{\mathop{\mathrm{Im}}}
\renewcommand{\i}{\mathrm{i}}

\newcommand{\CC}{{\mathbb C}}

\newcommand{\bM}{{\bf M}}
\newcommand{\bN}{{\bf N}}
\newcommand{\bI}{{\bf I}}

\begin{document}

\title{A boundary integral equation with the generalized Neumann kernel for the Ahlfors map}

\author{Mohamed M.S. Nasser$^{\rm a,b}$ and Ali H.M. Murid$^{\rm c,d}$}

\date{}
\maketitle
\vskip-0.7cm
\centerline{$^{\rm a}$Department of Mathematics, Faculty of Science, King Khalid University,} %
\centerline{P. O. Box 9004, Abha 61413, Saudi Arabia.}%
\centerline{$^{\rm b}$Department of Mathematics, Faculty of Science, Ibb University,} %
\centerline{P. O. Box 70270, Ibb, Yemen.} %
\centerline{E-mail: mms\_nasser@hotmail.com}

\vskip0.2cm
\centerline{$^{c}$Department of Mathematical Sciences, Faculty of Science, } %
\centerline{Universiti Teknologi Malaysia, 81310 UTM Johor Bahru, Johor, Malaysia.} %
\centerline{$^{d}$UTM Centre for Industrial and Applied Mathematics (UTM-CIAM), } %
\centerline{Universiti Teknologi Malaysia, 81310 UTM Johor Bahru, Johor, Malaysia.} %
\centerline{E-mail: alihassan@utm.my}

%\vskip1.2cm
\begin{center}
\begin{quotation}
{\noindent {\bf Abstract.\;\;}%
This paper presents a uniquely solvable boundary integral equation with the generalized Neumann kernel for the Ahlfors map of bounded multiply connected regions.  
}%
\end{quotation}
\end{center}
\begin{center}
\begin{quotation}
\end{quotation}
\end{center}

\begin{center}
\begin{quotation}
\end{quotation}
\end{center}
%

%--------------------------------------------------------------
\section{Introduction}
\label{sc:int}
%--------------------------------------------------------------

A boundary integral equation with the generalized Neumann kernel has been presented in~\cite{nas-fun,nas-siam,nas-jmaa,nas-jmaa2,nas-siam2} for computing the conformal mapping from multiply connected regions onto Koebe's $39$ canonical slit regions as well as the canonical region obtained by removing rectilinear slits from a strip. Only the right-hand side of the integral equation is different from one canonical region to another. In this paper, we shall extend the method presented in~\cite{nas-fun,nas-siam,nas-jmaa,nas-jmaa2,nas-siam2} to compute the Ahlfors map of bounded multiply connected regions. We shall show that the same integral equation used in~\cite{nas-fun,nas-siam,nas-jmaa,nas-jmaa2,nas-siam2} can be used to compute the Ahlfors map. However, the right-hand side of the integral equation contains the zeros of the Ahlfors map which are unknowns. 

A well known integral equation for the Ahlfors map is the boundary integral equation derived by Kerzman and Stein~\cite{ker78} for the Szeg\"o kernel. This integral equation has been used in~\cite{ker86,lee,odo} to compute the conformal mapping from simply connected regions onto the unit disc. A generalization of Kerzman-Stein method to compute the Ahlfors map for bounded multiply connected regions is given in Bell~\cite{bel86}. The method presented in~\cite{bel86} can be used to compute the Ahlfors map without relying on the zeros of the Ahlfors map. See also~\cite{mur,teg98,teg}. 

The approach proposed in this paper is of theoretical and practical interest. Firstly, it will be shown that the boundary integral equation with the generalized Neumann kernel used to compute the conformal mapping in~\cite{nas-fun,nas-siam,nas-jmaa,nas-jmaa2,nas-siam2} is valid for the Ahlfors map. Secondly, the approach proposed in this paper can be combined with the approach proposed in~\cite{bel86} to obtain a numerical method for computing the zeros of the Ahlfors map which is an interesting problem. The problem of computing the zeros of the Ahlfors map is more difficult than computing the Ahlfors map itself. A method for computing the zeros of the Ahlfors map for an annulus is presented in~\cite{teg98,teg}.

For other integral equations for the Ahlfors map of doubly connected regions, see~\cite{mur}. However, the right-hand side of these integral equations contains also an unknown constant which is a zero of the Ahlfors map.

%------------------------------------------------------------------
\section{Auxiliary material}
\label{sc:aux}
%------------------------------------------------------------------

Let $G$ be a bounded multiply connected region in the complex plane $\CC$ of connectivity $m\ge1$ with the boundary $\Gamma=\partial G$. The boundary $\Gamma$ consists of $m$ closed smooth Jordan curves $\Gamma_1, \Gamma_2, \ldots, \Gamma_m$ of which $\Gamma_m$ contains the others. The orientation of $\Gamma$ is such that $G$ is always on the left of $\Gamma$. For $j=1,\ldots,m$, the curve $\Gamma_j$ is parametrized by a $2\pi$-periodic twice continuously differentiable complex function $\eta_j(t)$ with non-vanishing first derivative
\[ 
\dot\eta_j(t)=d\eta_j(t)/dt\ne 0,\quad t\in J_j=[0,2\pi].
\]
The total parameter domain $J$ is the disjoint union of the intervals $J_1,\ldots,J_m$. We define a parametrization of the whole boundary $\Gamma$ as the complex function $\eta$ defined on $J$ by $\eta(t)=\eta_j(t)$ if $t\in J_j$, i.e.,
\begin{equation}\label{e:eta}
\eta(t)= \left\{ \begin{array}{l@{\hspace{0.5cm}}l}
\eta_{1}(t),&t\in J_{1}=[0,2\pi],\\
\hspace{0.3cm}\vdots\\
\eta_m(t),&t\in J_m=[0,2\pi].
\end{array}
\right .
\end{equation}
The meaning of the notation in~(\ref{e:eta}) is as follows: For a given $\hat t\in[0,2\pi]$, to evaluate the value of $\eta(t)$ at $\hat t$, we should know in advance the interval $J_j$ to which $\hat t$ belongs, i.e., we should know the boundary $\Gamma_j$ contains $\eta(\hat t)$, then we compute $\eta(\hat t)=\eta_j(\hat t)$.

Let $H$ be the space of all real H\"older continuous functions on $\Gamma$. In view of the smoothness of $\eta$, a real H\"older continuous function $\hat\phi$ on $\Gamma$ can be interpreted via $\phi(t)=\hat\phi(\eta(t))$, $t\in J$, as a function $\phi\in H$; and vice versa. Let $S$ be the subspace of $H$ that consists of real piecewise constant functions of the form
\begin{equation}\label{e:piece}
h(t) = \left\{
\begin{array}{l@{\hspace{0.5cm}}l}
 h_1,     &t\in J_1, \\
 \vdots   & \\
 h_m,     &t\in J_m, \\
\end{array}%
\right.
\end{equation}
with real constants $h_1,\ldots,h_m$. 

For a fixed point $a\in G$, we define a complex-valued function $A(t)$ on $\Gamma$ by
\begin{equation}\label{e:A}
A(t) = \eta(t)-a,
\end{equation}
The generalized Neumann kernel formed with $A$ is defined by
\begin{equation}\label{e:N}
 N(s,t) =  \frac{1}{\pi}\Im\left(
 \frac{A(s)}{A(t)}\frac{\dot\eta(t)}{\eta(t)-\eta(s)}\right).
\end{equation}
We define also a real kernel $M$ by
\begin{equation}\label{e:M}
 M(s,t) =  \frac{1}{\pi}\Re\left(
 \frac{A(s)}{A(t)}\frac{\dot\eta(t)}{\eta(t)-\eta(s)}\right).
\end{equation}
The kernel $N$ is continuous and the kernel $M$ has a cotangent singularity type. Hence, the operator
\begin{equation}\label{e:opN}
  \bN \mu(s) = \int_J N(s,t) \mu(t) dt, \quad s\in J,
\end{equation}
is a Fredholm integral operator and the operator
\begin{equation}\label{e:opM}
  \bM\mu(s) = \int_J  M(s,t) \mu(t) dt, \quad s\in J,
\end{equation}
is a singular integral operator. For more details, see~\cite{mur03,wegm}.

%-------------------------------------------------------------
\section{The Ahlfors Map}
\label{sc:map}
%-------------------------------------------------------------

The Ahlfors map $\omega$ is a branched $m$-to-one analytic function mapping $G$ onto the unit disk. The Ahlfors map takes each boundary component $\Gamma_j$ one-to-one and onto the unit circle. With the normalization 
\[
\omega(a)=0, \quad \omega'(a)>0,
\]
the Ahlfors map is unique. For more details, see~\cite[p.~94]{kra} and~\cite[p.~379]{neh}.

The Ahlfors map has exactly $m$ zeros in $G$. One of these zeros is $a$ and the remaining $m-1$ zeros are unknowns. We represent these unknown zeros as $a_1,a_2,\ldots,a_{m-1}$. For $j=1,2,\ldots,m-1$, suppose that $z_j$ is a fixed point inside $\Gamma_j$. Hence, the Ahlfors map $\omega$ can be written as
\begin{equation}\label{e:om}
\omega(z) = c(z-a)\prod_{j=1}^{m-1}\frac{a-z_j}{a-a_j}\,\prod_{j=1}^{m-1}\frac{z-a_j}{z-z_j}\,e^{(z-a)f(z)},
\end{equation}
where $c=\omega'(a)>0$ and $f$ is an auxiliary analytic function in $G$. Introducing the product $\prod_{j=1}^{m-1}\frac{z-a_j}{z-z_j}$ in~(\ref{e:om}) is necessary to ensure that both sides of~(\ref{e:om}) have the same winding numbers along the boundary components since this is important for the existence of an analytic function $f(z)$ such that~(\ref{e:om}) is valid. To ensure that $\omega'(a)=c$, we introduce the product $\prod_{j=1}^{m-1}\frac{a-z_j}{a-a_j}$ in~(\ref{e:om}).  

Since $|\omega(\eta(t))|=1$ on the boundary $\Gamma$, we obtain
\[
1=c|\eta(t)-a|\prod_{j=1}^{m-1}\frac{|a-z_j|}{|a-a_j|}\,\prod_{j=1}^{m-1}\frac{|\eta(t)-a_j|}{|\eta(t)-z_j|}\,e^{\Re[A(t)f(\eta(t))]}
\]
which implies that the boundary values of the function $f$ satisfy
\begin{equation}\label{e:RH-Ah}
\Re[A(t)f(\eta(t))]= -\ln c-\sum_{j=1}^{m-1}\ln\frac{|a-z_j|}{|a-a_j|}-\ln|\eta(t)-a|
-\sum_{j=1}^{m-1}\ln|\eta(t)-a_j|+\sum_{j=1}^{m-1}\ln|\eta(t)-z_j|.
\end{equation}
If the zeros $a_1,a_2,\ldots,a_{m-1}$ are known, then the function $f$ is the unique solution of the RH problem
\begin{equation}\label{e:RH-Ah2}
\Re[A(t)f(\eta(t))]= h(t)+\gamma(t)
\end{equation}
where
\begin{equation}\label{e:RH-gam}
\gamma(t)=-\ln|\eta(t)-a|-\sum_{j=1}^{m-1}\ln|\eta(t)-a_j|+\sum_{j=1}^{m-1}\ln|\eta(t)-z_j|
\end{equation}
and
\begin{equation}\label{e:RH-h}
h(t)=  -\ln c-\sum_{j=1}^{m-1}\ln\frac{|a-z_j|}{|a-a_j|}\in S.
\end{equation}

Let $\mu(t)=\Im[A(t)f(\eta(t))]$, i.e., the boundary values of the function $f(z)$ are given by 
\begin{equation}\label{e:Af}
Af=\gamma+h+\i\mu
\end{equation}
where $\mu$ and $h$ are unknown functions. It follows from~\cite{nas-fun} that the function $\mu$ is the unique solution of the integral equation
\begin{equation}\label{e:ie}
(\bI-\bN)\mu=-\bM\gamma
\end{equation}
and the function $h$ is given by
\begin{equation}\label{e:h}
h=[\bM\mu-(\bI-\bN)\gamma]/2.
\end{equation}
By obtaining the functions $\mu$ and $h$, we obtain the boundary values of the function $f(z)$ from~(\ref{e:Af}). Then, the values of the function $f(z)$ can be computed for $z\in G$ by the Cauchy integral formula
\begin{equation}\label{e:f-cau}
f(z)=\frac{1}{2\pi\i}\int_{\Gamma} \frac{\gamma+h+\i\mu}{A} \frac{1}{\eta-z}d\eta.
\end{equation}
Then, the Ahlfors map $\omega(z)$ is computed from~(\ref{e:om}).

%-------------------------------------------------------------
\section{Numerical Example}
\label{sc:exm}
%-------------------------------------------------------------

A reliable procedure for solving the integral equation~(\ref{e:ie}) numerically is by using the Nystr\"om method with the trapezoidal rule with $n$ equidistant nodes in each interval $J_j$, $j=1,\ldots,m$. See~\cite{nas-fun,nas-siam,nas-jmaa,nas-jmaa2,nas-siam2} for more details.  

\begin{figure}[ht] %
\centerline{
\scalebox{0.29}{\includegraphics{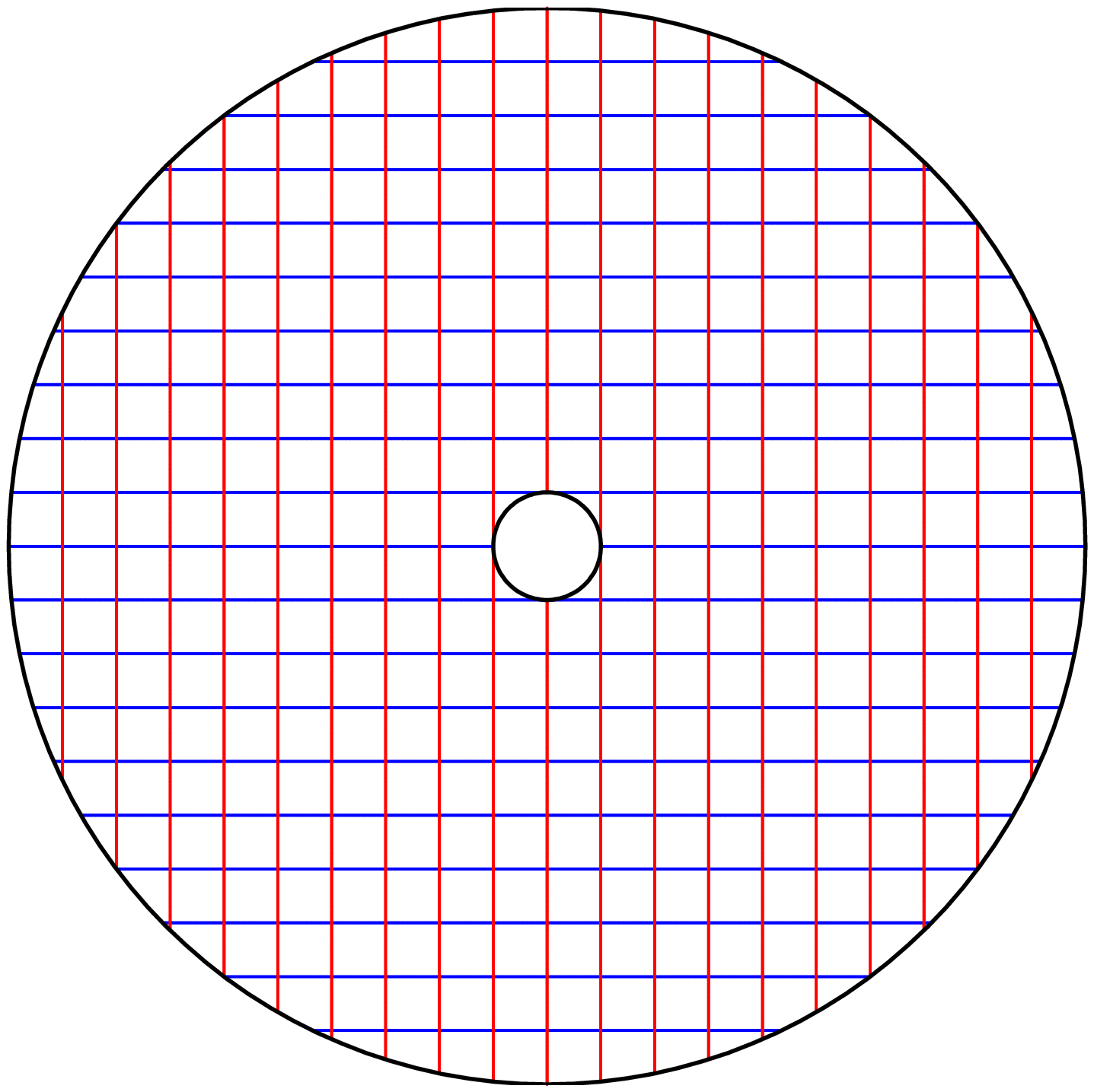}}\hfill
\scalebox{0.29}{\includegraphics{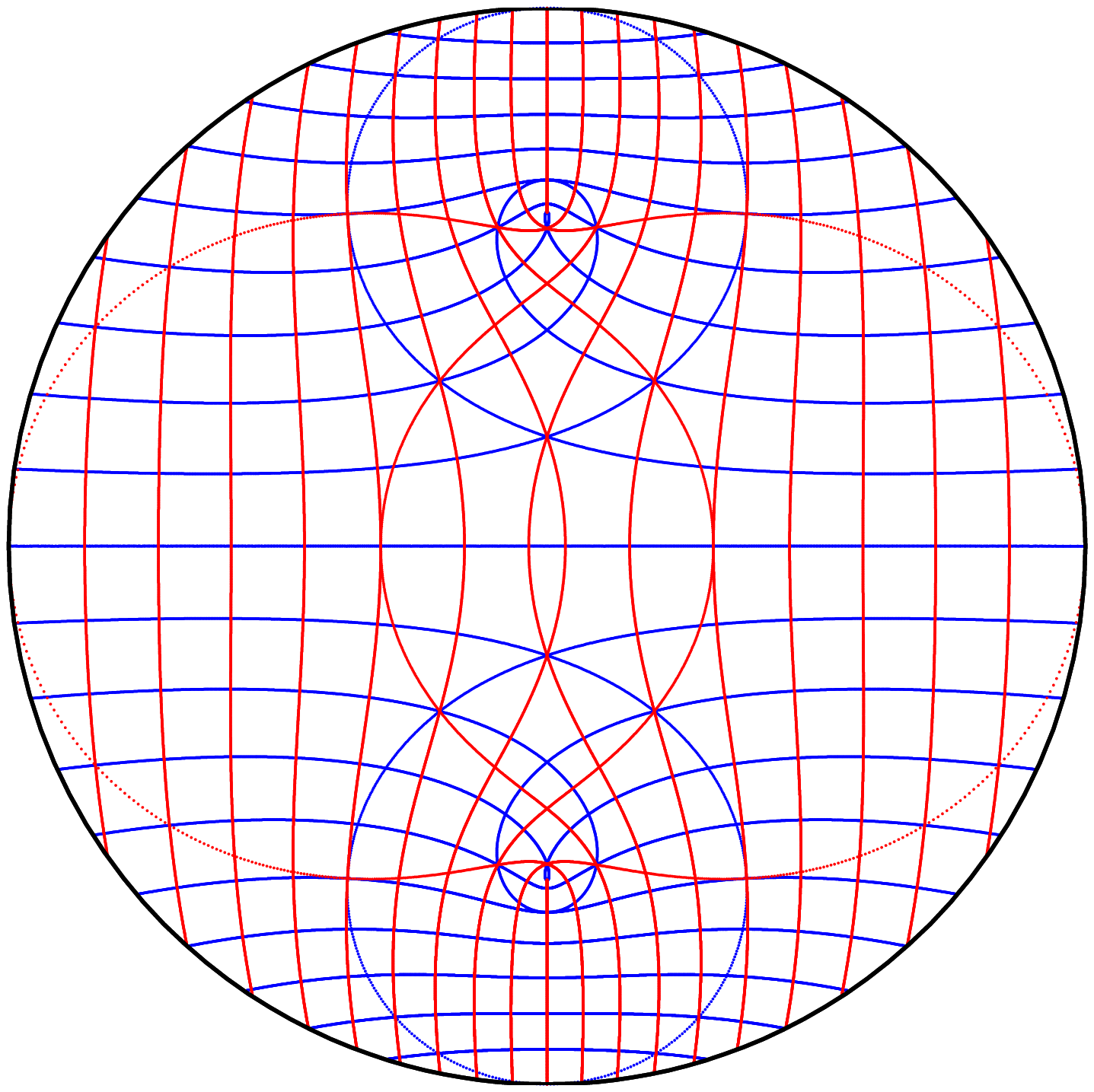}}\hfill
\scalebox{0.29}{\includegraphics{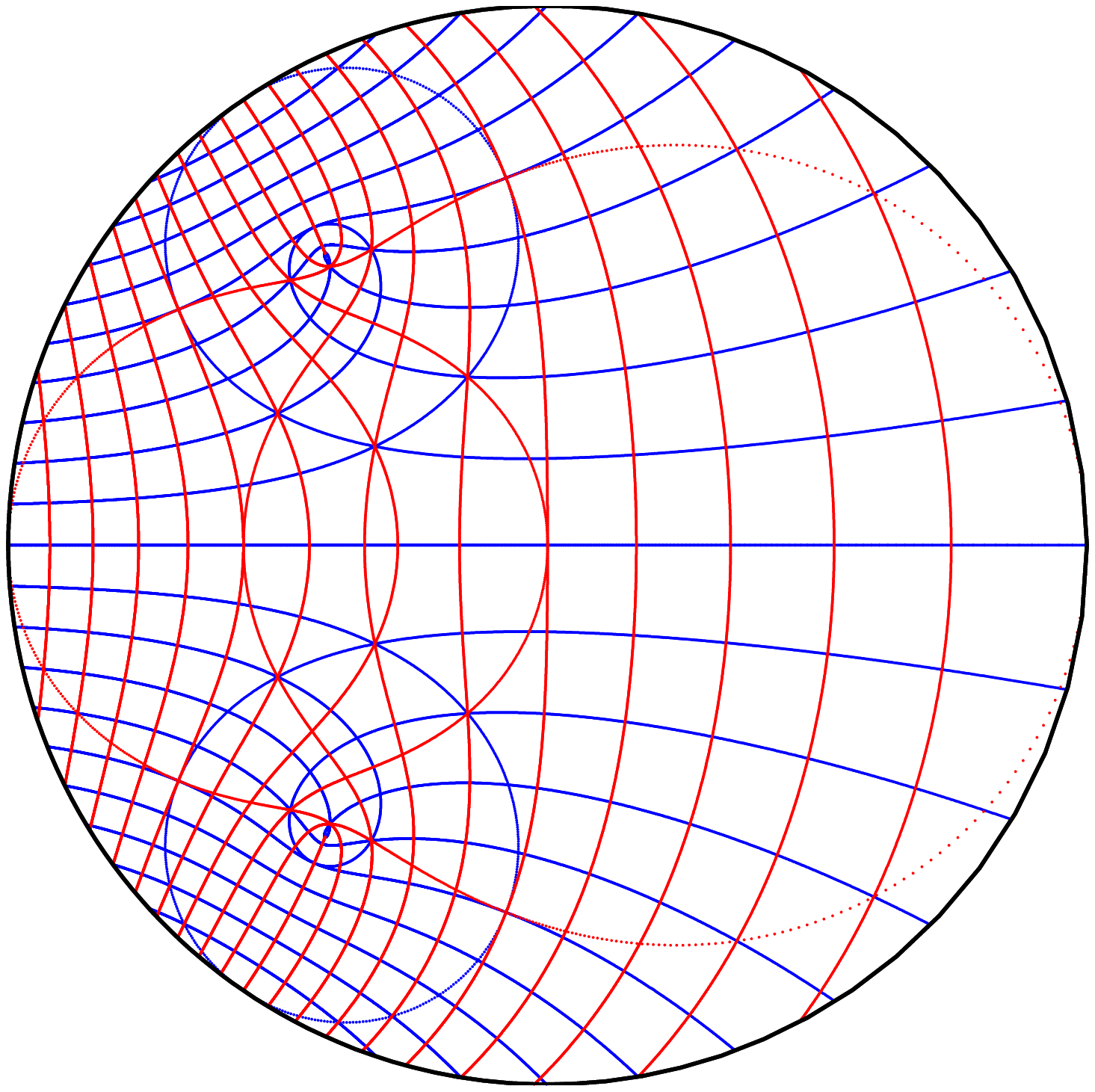}}
}
\centerline{\hfill (a) \hfill\hfill (b) \hfill\hfill (c)\hfill}
\caption{(a) Grid over the original region $G$. (b) The image of the region $G$ for $a=\frac{1}{\sqrt{10}}$. (c) The image of the region $G$ for $a=0.5$.} 
\label{f:exm}
\end{figure}

In this section, we shall consider a numerical example in the annulus $0.1<|z|<1$ ($m=2$). This example have been considered in~\cite{teg} where the Ahlfors map was computed using the Szeg\"o and the Garabedian kernels. The first zero of the Ahlfors map, $a$, is assumed to be an arbitrary point in $G$. The second zero of the Ahlfors map, $a_1$, was computed in~\cite{teg} as the zero of the Szeg\"o kernel. It was proved in~\cite[Theorem~2]{teg} that $a_1=-0.1/\overline{a}$. Orthogonal grid over the original region $G$ is shown in Fig.~\ref{f:exm}(a). The image of the region $G$ is shown in Fig.~\ref{f:exm}(b) for $a=\frac{1}{\sqrt{10}}$ (hence $a_1=\frac{-1}{\sqrt{10}}$) and in Fig.~\ref{f:exm}(c) for $a=0.5$ (hence $a_1=-0.2$). The same number, $n=128$, of nodes in each interval $J_k$, $k=1,2$, was used.

%--------------------------------------------------------------------

\end{document}